\documentclass{article}
\usepackage[utf8]{inputenc}
\usepackage{amsthm}
\usepackage{amsmath}
\usepackage{amssymb}
\usepackage{pgfplots}
\pgfplotsset{width=10cm,compat=1.9}
\usepgfplotslibrary{external}
\usepackage{geometry}
\usepackage{listings}
\usepackage{graphicx}
\usepackage{subcaption}
\usepackage{wrapfig,lipsum,booktabs}
\newtheorem{theorem}{Theorem}[section]

\usepackage[
backend=biber,
style=alphabetic,
sorting=ynt
]{biblatex}
\usepackage{graphicx} % Required for inserting images

\title{A Family of New Formulas for the Euler-Mascheroni Constant}
\author{Noah Ripke}
\date{}

\begin{document}

\maketitle
\begin{abstract}
We introduce and prove several new formulas for the Euler-Mascheroni Constant. This is done through the introduction of the defined E-Harmonic function, whose properties, in this paper, lead to two novel formulas, alongside a family of formulas. While the paper does introduce many new approximations, it does not exhaust the possibilities of the E-Harmonic function but provides a strong first dive into its natural conclusions. We hope that the diversity of new formulas may provide stepping stones to a proof (or disproof) of the irrationality of the Euler-Mascheroni constant.\newline \newline \newline \newline
\end{abstract}
\section{Introduction}
\textrm{The Euler-Mascheroni constant is the limiting difference between the harmonic series and the natural logarithm \cite{1} proving a fascinating connection between the discrete and continuous. }
\[
 \lim_{n \rightarrow \infty}\sum_{i=1}^{n}{\frac{1}{i}} - \ln{n} = \gamma
\]

\textrm{But the constant's claim to fame is that a proof (or disproof) of its irrationality has eluded Mathematicians for centuries. Several formulas have been discovered in the hope of finding a path towards irrationality, and faster convergence. This paper will introduce two new concise and effective formulas that may be used towards a theory of irrationality, and the methods used to discover these formulas can be replicated to create an infinite family of formulas (Most of which are slightly more convoluted). The star of this paper is the E-Harmonic function \(E(n)\) whose interesting properties will hopefully pave a path towards many new discoveries involving the Euler-Mascheroni constant.}
\newpage
\section{The E-Harmonic Function}
\textrm{The limiting connection between the Harmonic Series and the natural logarithm leads us to the intuitive analysis of what I call the E-Harmonic function: }

\begingroup
\Large
\[
E(n) = e^{H_n}
\]
\endgroup

\textrm{Note that \(H_n\) is simply \(\sum_{i=1}^{n}{\frac{1}{i}}\). The best way to understand the power of the E-Harmonic Function is to see its graph:}

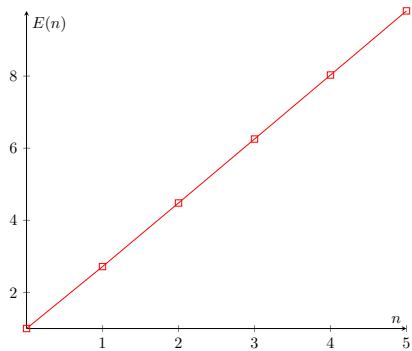
\begin{center}
\begin{tikzpicture}[scale=0.6]
\centering
\begin{axis}[
    axis lines = center,
    xlabel = \(n\),
    ylabel = {\(E(n)\)},
]

\addplot[
    color=red,
    mark=square,
    ]
    coordinates {
    (0,1)(1,2.718)(2, 4.482)(3,6.255)(4,8.031)(5,9.809)
    };
\end{axis}
\end{tikzpicture}
\end{center}

\textrm{The linear-looking aspect of this graph should be somewhat expected due to the Euler-Mascheroni Constant, but the rate at which the E-Harmonic function approaches linear is what makes the formulas derived from it converge so quickly. When first looking at this graph, the next intuitive step is to ask what the limiting difference between each successive term is: and it happens to be \(e^{\gamma}\)!}

\begin{theorem}[E-Harmonic Linearization Theorem]
An integral theorem to all E-Harmonic derived formulas for the Euler-Mascheroni constant
\end{theorem}

\begingroup
\Large
\[
\lim_{n \rightarrow \infty}E(n) = \lim_{n \rightarrow \infty}e^{\gamma}n
\]
\endgroup

\begin{proof}
A simple proof comes from the original definition for the Euler-Mascheroni constant:
\[
 \lim_{n \rightarrow \infty}H_n - \ln{n} = \gamma
\]
By raising each side to the power of \(e\), we receive this true statement:
\[
 \lim_{n \rightarrow \infty}\frac{E(n)}{n} = e^\gamma
\]
Trivial application of limit theory will obtain the desired result.
\end{proof}

\textrm{The power of the E-Harmonic Linearization Theorem is our ability to now substitute the limiting behavior of \(E(n)\) with a simplistic function involving \(\gamma\). This will prove to be endlessly fruitful in the following sections.}

\section{Two Novel Formulas}
\textrm{The following two sub-sections will show the two main results of this paper. These formulas, derived from \(E(n)\), are both interesting in their own right and are generally the most simplistic of all the other formulas described in this paper (And thus will be discussed separately).}
\subsection{\(E(n+1)-E(n)\)}
\textrm{In this subsection, we will prove and discuss the following formula:}
\begin{equation}
\large
\lim_{n \rightarrow \infty}H_n + \ln{(e^{\frac{1}{n+1}}-1)} = \gamma
\end{equation}
\begin{proof}
Using Theorem 2.1, we start with the following true statement:
\[
\lim_{n \rightarrow \infty}E(n+1)-E(n) = e^\gamma
\]
\textrm{By rewriting \(E(n)\) as \(e^{H_n}\), and factoring out a common \(E(n)\): }
\[
\lim_{n \rightarrow \infty}e^{H_n}(e^{\frac{1}{n+1}}-1) = e^\gamma
\]
\textrm{Finally, when we take the natural log of both sides, we obtain the desired result.}
\end{proof}
\textrm{This formula is interesting in a few aspects. For one, its similarity to the original definition of the Euler-Mascheroni constant allows us to notice that the functions \(f(x) = e^{\frac{1}{x+1}}-1\), and \(g(x) = \frac{1}{x}\) behave very similarly as \(x\) approaches infinity. This might be obvious from an early look at the Taylor Series expansion for \(f(x)\), but it leads us to an interesting observation.\newline}

\textrm{In a philosophical manner, it seems the original definition of the Euler-Mascheroni constant is a simple byproduct (and approximation) of formula 3.1.1. While formula 3.1.1 is more complex than the original simplistic nature of \(\gamma\), it seems to be the "father" of our primary definition. And in fact, this new formula converges much faster! To use in actual computations, one can create an \(n\)-degree Taylor Series for \(e^{\frac{1}{x+1}}-1\), and choose some large \(n\); the larger the \(n\), the faster the convergence. Of course, this approximation would also involve a second Taylor Series for \(\ln{(x)}\) for small \(x\) approaching 0.\newline}

\textrm{Hopefully, this formula may also be used as a device towards a theory of irrationality of \(\gamma\). For one, this formula separates itself from other equations in that it involves \(e\), and not some value shooting off towards infinity, but a function of \(e\) that converges to 0.}
\newpage
\subsection{\(E(n+1)*E(n)\)}
\textrm{In this subsection, we will prove and discuss the following formula:}
\begin{equation}
\large
\lim_{n \rightarrow \infty}H_n + \frac{1}{2(n+1)} - \frac{1}{2}\ln{(n^2+n)} = \gamma
\end{equation}
\begin{proof}
\textrm{In a similar manner to the proof in 3.1, we start by rewriting \(E(n+1)E(n)\) using Theorem 2.1:}
\[
\lim_{n \rightarrow \infty}E(n+1)E(n) = \lim_{n \rightarrow \infty} e^{2\gamma}(n^2+n)
\]
\textrm{Since \(H_{n+1} = H_n + \frac{1}{n+1}\), we can take the natural log of both sides to obtain this statement:}
\[
\lim_{n \rightarrow \infty}2H_n + \frac{1}{n+1} - \ln{(n^2+n)} = 2\gamma
\]
\textrm{The equation is perfectly fine in this format, but by dividing by two, we can reach the desired result.}
\end{proof}

\textrm{While this formula seems interesting, it actually converges just as fast as the original definition for \(\gamma\). This may be disappointing until an observation is made. We seem to have a "leftover" term of \(\frac{1}{2(n+1)}\), which converges to 0 of course as \(n\) approaches infinity. It happens to be that removing this term allows for significantly faster convergence. \newline \newline}

\begin{center}
\begin{tabular}{|c||c|c||} 

 \hline
\(n\) & With "Leftover" Term & Without "Leftover" Term\\ [0.5ex] 
 \hline\hline
 1 & 0.90342640972 & 0.65342640972\\
 \hline
 10 & 0.624182616527 & 0.578728071072\\ 
 \hline
 100 & 0.582182661274 & 0.577232166225\\ 
 \hline
 1,000 & 0.577715331901 & 0.577215831402\\ 
 \hline
 10,000 & 0.577265661568 & 0.577215666568\\ 
 \hline
\end{tabular}

\end{center}
\textrm{\newline \newline}
\textrm{This idea doesn't come from nowhere. If you were to combine our original formula for \(\gamma\) and a first-degree approximation of the formula derived in 3.1, you would achieve the same equation! This equation will be expanded upon in the next section which aims to discuss these formulas in a more broad, and generalized way.}
\newpage
\section{A Family of Formulas}
We generalize the formulas found in the previous sections and notice interesting results that can be derived from them.

\subsection{\(E(n+k)-E(n)\)}

This is a clear continuation of the idea brought up in 3.1, and gives rise to the following equation:

\begin{equation}
\large
\lim_{n \rightarrow \infty} H_n+ \ln(e^{\sum_{i=1}^{k}\frac{1}{n+i}}-1) -\ln(k) = \gamma
\end{equation}
\begin{proof}
\textrm{Clearly, by applying Theorem 2.1:}
\[
\lim_{n \rightarrow \infty}E(n+k)-E(n) = ke^{\gamma}
\]
\textrm{Thus, by factoring the left side into:}
\[
\lim_{n \rightarrow \infty}E(n)(e^{\frac{1}{n+1}+...+\frac{1}{n+k}}-1) = ke^{\gamma}
\]
\textrm{We see that our desired formula arises easily by taking the natural log of both sides. We will call this equation: }
\[
L(n,k)=H_n+ \ln(e^{\sum_{i=1}^{k}\frac{1}{n+i}}-1) -\ln(k)
\]
\end{proof}

Now, clearly from construction, as \(n\) limits to \(\infty\) we get \(\gamma\), but is the same true for constant \(n\) and limiting \(k\)? We verify:

\[
\lim_{k \rightarrow \infty} L(n,k) = \gamma
\]
\begin{proof}
\textrm{By working backwards from the definition of \(L(n,k)\), we let:}
\[
\lim_{k \rightarrow \infty}H_n+ \ln(e^{\sum_{i=1}^{k}\frac{1}{n+i}}-1) -\ln(k) = L
\]
\textrm{Thus by raising both sides to the \(e\):}
\[
\lim_{k \rightarrow \infty}\frac{E(n+k)-E(n)}{k} = e^L
\]
\textrm{Since \(E(n+k)\) limits to infinity, we can apply Theorem 2.1 to see that:}
\[
\lim_{k \rightarrow \infty}\frac{e^{\gamma}(n+k)-E(n)}{k} = e^L
\]
\textrm{Then after basic re-arranging we find that: }
\[
e^L = e^{\gamma} - \lim_{k \rightarrow \infty}\frac{e^{\gamma}n-E(n)}{k}
\]
\textrm{Since our lasting limit clearly goes to 0 since the numerator is a constant, we find that \(L=\gamma\).}
\end{proof}

\textrm{Now from the proof, it's evident that for larger \(n\), our approximation can become more accurate since \(|e^{\gamma}n-E(n)|\) is decreasing. But one can conveniently choose \(n=0\) and out comes a completely new formula:}

\begin{equation}
\large
\lim_{k \rightarrow \infty} \ln(E(k)-1)-\ln(k) = \gamma
\end{equation}

\textrm{The formula is amazingly close to the original definition of \(\gamma\), and yet, can get away with being much more accurate.}

\subsection{\(E(n)*E(n+1)*...*E(n+k)\)}
\textrm{There are multiple natural continuations of the product formula, but this one I find to be most interesting.}

\begin{equation}
\large
\lim_{n \rightarrow \infty} H_{n+k} - \frac{1}{k+1}\sum_{i=1}^{k}{\frac{i}{n+i}}-\frac{1}{k+1}\ln(\frac{(n+k)!}{(n-1)!}) = \gamma
\end{equation}

\begin{proof}
\textrm{We start out with our product and note that using Theorem 2.1: }
\[
\lim_{n \rightarrow \infty}E(n)*...*E(n+k)=\lim_{n \rightarrow \infty}e^{(k+1)\gamma}n(n+1)...(n+k)
\]
\textrm{Then by taking the natural log of both sides, we see that after grouping the \((k+1)H_n\) and the leftover fractional terms into a sum, we get: }
\[
\lim_{n \rightarrow \infty} H_n + \frac{1}{k+1}\sum_{i=1}^{k}{\frac{k-i+1}{n+i}}-\frac{1}{k+1}\ln(\prod_{i=0}^{k}(n+i)) = \gamma
\]
\textrm{Then finally, turning our product into a factorial, and re-arranging our sum, we get the desired result.}
\end{proof}

\textrm{For fixed \(k\), we can see that our approximation will approach \(\gamma\), but as was discussed in the previous subsection, does increasing \(k\) also limit to \(\gamma\)? Calling our new equation \(P(n,k)\), we verify:}

\[
\large
\lim_{k \rightarrow \infty} P(n,k) = \gamma
\]
\begin{proof}
\textrm{We start with our limit, and let it equal \(L\), from there we will prove that \(L=\gamma\)}
\[
\lim_{k \rightarrow \infty} H_{n+k} - \frac{1}{k+1}\sum_{i=1}^{k}{\frac{i}{n+i}}-\frac{1}{k+1}\ln(\frac{(n+k)!}{(n-1)!}) = L
\]
\textrm{Rewrite the equation as such: }
\[
\lim_{k \rightarrow \infty} H_{n+k} - \frac{k}{k+1}+\frac{1}{k+1}\sum_{i=1}^{k}{\frac{n}{n+i}}-\frac{1}{k+1}\ln((n+k)!)+\frac{1}{k+1}\ln((n-1)!) = L
\]
\textrm{Since \((n-1)!\) is constant, and our sum will clearly go to 0, we get: }
\[
\lim_{k \rightarrow \infty} H_{n+k} - 1 -\frac{1}{k+1}\ln((n+k)!) = L
\]
\textrm{Now, we use Stirling's Approximation for \(\ln(n!) \approx n\ln(n)-n\). This is rigorous since \((n+k)\) approaches infinity, so as long as we maintain our limit, our statement still holds.}
\[
\lim_{k \rightarrow \infty} H_{n+k} - 1 -\frac{1}{k+1}((n+k)\ln(n+k)-(n+k)) = L
\]
\textrm{Since \(\frac{n+k}{k+1}\) approaches 1, we get:} 
\[
\lim_{k \rightarrow \infty} H_{n+k} -\ln(n+k) = L
\]
\textrm{But this is simply the limiting difference between the Harmonic Series and the natural log, which is the original definition of \(\gamma\)! Thus \(L=\gamma\).}
\end{proof}

\textrm{Now that we have this powerful fact, we can use \(n=1\) as our convenient value, and quickly derive this new formula:}
\[
\large
\lim_{k \rightarrow \infty} H_{k+1} + \frac{H_{k+1}}{k+1} - 1 - \frac{1}{k+1}\ln((k+1)!) = \gamma
\]
\textrm{But in this case, it may be more easily understood to replace \(k+1\) with simply \(k\) to get:}
\begin{equation}
\large
\lim_{k \rightarrow \infty} H_{k} + \frac{H_{k}}{k} - 1 - \frac{1}{k}\ln(k!) = \gamma
\end{equation}
\textrm{Now, one can apply Stirling's approximation to our lasting factorial term, but the resulting formula is one that is not too interesting, nor helpful for approximations. \newline}

\textrm{As can be seen, the formulas shown in this paper just scratch the surface of the possibilities for new formulas. Abusing properties of linearity and products are just natural first steps, but further research can be done on new approaches to the E-Harmonic function that may provide even more insightful results. Such approaches may include unique sums, further study of \(L(n,k)\) and \(P(n,k)\), etc.}

\printbibliography

@article{1,
    author = {L. Euler},
    title = {De progressionibus harmonicis observationes},
    journal = {Commentarii academiae scientiarum Petropolitanae},
    year = {1740},
    volume = {7},
    number = {},
    pages = {150-161},
    doi = {},
    url = {},
}

\end{document}